\documentclass[reqno]{amsart}
\usepackage{amssymb}

\newcommand{\zero}{\mathbf{0}}
\newcommand{\one}{\mathbf{1}}

\newcommand{\V}{\mathcal{V}}

\newcommand{\C}{\mathbf{C}}
\newcommand{\B}{\mathbf{B}}
\newcommand{\E}{\mathbf{E}}
\newcommand{\F}{\mathbf{F}}

\newcommand{\ba}{\boldsymbol{a}}
\newcommand{\bb}{\boldsymbol{b}}
\newcommand{\bc}{\boldsymbol{c}}
\newcommand{\bd}{\boldsymbol{d}}
\newcommand{\be}{\boldsymbol{e}}

\newcommand{\sd}{\smallsetminus}
\newcommand{\ch}{\mathfrak{C}_2}
\numberwithin{equation}{section}

\theoremstyle{plain}%default
\newtheorem{theorem}{Theorem}[section]
\newtheorem{proposition}[theorem]{Proposition}
\newtheorem{lemma}[theorem]{Lemma}
\newtheorem{corollary}[theorem]{Corollary}

\theoremstyle{definition}

\newtheorem{definition}[theorem]{Definition}

\theoremstyle{remark}
\newtheorem*{remark}{Remark}

\begin{document}

\title[Non-distributive varieties]%
{Congruence lattices of free lattices in non-distributive
varieties}
 \author[M.~Plo\v s\v cica]{Miroslav~Plo\v s\v cica}
 \thanks{The research of the first author was
         supported by Slovak GAV grant 1230/96.}
 \address{Mathematical Institute\\
         Slovak Academy of Sciences\\
         Gre\v s\'akova 6\\
         04001 Ko\v sice\\
         Slovakia} \email{ploscica@linux1.saske.sk}

 \author[J.~T\accent23uma]{Ji\v r\'\i~T\accent23uma}
 \thanks{The second author was partially supported by GA\v CR
         grant no. 201/95/0632 and by GA AV \v CR grant no. A
         1019508.}
 \address{Department of Algebra\\
          Faculty of Mathematics and Physics\\
          Sokolovsk\'a 83\\
          Charles University\\
          186 00 Praha 8\\
          Czech Republic}
 \email{tuma@karlin.mff.cuni.cz}

\author[F.~Wehrung]{Friedrich~Wehrung}
\address{C.N.R.S.\\
         D\'epartement de Math\'ematiques\\
         Universit\'e de Caen\\
         14032 Caen Cedex\\
         France}
 \email{gremlin@math.unicaen.fr}

 \date{\today}

\begin{abstract}
We prove that for any free lattice $F$ with at least
$\aleph_2$ generators in any non-distributive variety of
lattices, there exists no sectionally complemented
lattice $L$ with congruence lattice isomorphic to the one
of $F$. This solves a question formulated by Gr\"atzer and
Schmidt in 1962. This yields in turn further examples of
simply constructed distributive semilattices that are not
isomorphic to the semilattice of finitely generated two-sided
ideals in any von Neumann regular ring.
\end{abstract}

\maketitle

\section*{Introduction}

One of the oldest and most famous open problems in lattice
theory, the \emph{Congruence Lattice Problem}, is to
decide whether for every distributive (join-) semilattice $S$
with zero, there exists a lattice $L$ such that the semilattice
\(\C(L)\) of compact congruences of $L$ (the
\emph{congruence semilattice} of $L$) is isomorphic to
$S$. Although the answer is not known yet, many partial results
have been obtained (see \cite{Tisc94} for a survey).
Among these are positive solutions of the Congruence
Lattice Problem in the case where $S$ has size at most
$\aleph_1$, or is a distributive lattice. In addition, it
turns out that in several cases, the solution lattice $L$ to
the problem is
\emph{sectionally complemented} (\emph{e.g.}, for the
finite case, see \cite{GrSc62}; the case where $S$ is
\emph{countable} results from unpublished work of G. M.
Bergman \cite{Berg} and results in \cite{Wehr2}). As there
seems to be a growing evidence that in all
\emph{known} cases, there exists a sectionally
complemented solution lattice $L$, one may be tempted to
formulate the even stronger conjecture that every
distributive semilattice with zero is isomorphic to the
congruence semilattice of a \emph{sectionally
complemented} lattice. This conjecture had in fact already
been formulated in \cite[Problem II.8]{Grat78}.

In \cite{Wehr2}, using a construction presented in
\cite{Wehr1}, the third author proves that it cannot be so,
by giving a distributive semilattice of size $\aleph_2$
that is not isomorphic to the congruence semilattice of
any lattice that is, in the terminology of \cite{Wehr2},
\emph{congruence splitting}. In particular, every lattice
which is either sectionally complemented, or relatively
complemented, or a direct limit of atomistic lattices, is
congruence splitting. Therefore, if one could prove that
for every lattice $L$ there exists a sectionally
complemented (or, more generally, congruence splitting)
lattice $K$ such that \(\C(L)\cong\C(K)\), then one would
obtain a negative solution to the Congruence Lattice Problem.

This turns out to be also an open problem, more specifically
the second part of \cite[Problem 1, p. 181]{GrSc62}. In this
paper we give a strong negative solution to this problem, by
proving (Corollary~\ref{C:FreeNon2}) that in any
non-distributive variety of lattices,
if $F$ is any (bounded or not) free lattice with at least
$\aleph_2$ generators, then there exists no congruence
splitting lattice $L$ such that \(\C(F)\cong\C(L)\); in
particular, $F$ has no con\-gru\-en\-ce-pre\-ser\-ving
embedding into any sectionally complemented lattice. By
earlier results in
\cite{Wehr2}, this implies that \(\C(F)\) is never
isomorphic to the semilattice of finitely generated two-sided
ideals in a von Neumann regular ring.
The restrictions on the lattice variety
are optimal, because of the classical result that says that
every distributive lattice embeds
con\-gru\-en\-ce-pre\-ser\-ving\-ly into a generalized Boolean
algebra.

The strategy of the proof is the following: by the results of
\cite{Wehr2}, the congruence semilattice of any congruence
splitting lattice satisfies a certain infinite axiom, the
\emph{Uniform Refinement Property} (URP). In this paper, we
introduce a slight weakening of URP, the \emph{weak Uniform
Refinement Property} (WURP), that is not satisfied by the
congruence semilattice of any free lattice with at least
$\aleph_2$ generators in any non distributive lattice variety
$\V$. The two cases in which the proof splits, namely whether
the diamond $M_3$ or the pentagon $N_5$ belongs to $\V$, are
treated in a similar fashion: they are decorated with three
$2$-element chains that somewhat concentrate into a finite
pattern the combinatorial core of the original infinite
WURP. As in \cite{Wehr1}, the reduction of the
infinite case to the finite case is done \emph{via}
Kuratowski's free set property (\cite{Kura51}; see also
\cite[Proposition 2.5]{Wehr1} for a short proof).

\section*{Notation and terminology}

We consider semilattices of compact congruences of lattices.
The semilattices are join semilattices with $0$. The mapping
assigning to every lattice $L$ its congruence semilattice
$\C(L)$ can be extended to a functor from the category of
lattices and lattice homomorphisms to the category of
semilattices with homomorphisms of
semilattices; in addition, this functor preserves direct
limits. The least and largest congruence on $L$ will be
respectively denoted by $\zero$ and $\one$.

For all elements
$a$ and $b$ of a lattice $L$, we will denote by $\Theta(a,b)$
the least congruence containing the pair $(a,b)$ and we will
then put \(\Theta^+(a,b)=\Theta(a\wedge b,a)\); thus
\(\Theta^+(a,b)\) is the least congruence $\theta$ on $L$
such that \(\theta(a)\leq\theta(b)\).

We say that a
homomorphism of semilattices \(\mu:\,S\to T\) is
\emph{weak-distributive} \cite[Section 1]{Wehr2} when
for all \(e\in S\) and all \(b_0,b_1\in T\) such that
\(\mu(e)=b_0\vee b_1\), there are elements \(a_0\) and \(a_1\)
of \(S\) such that \(\mu(a_0)\leq b_0\), \(\mu(a_1)\leq b_1\)
and \(e=a_0\vee a_1\).

For every non negative integer $n$, we
will identify $n$ with the finite set
(initial ordinal) \(\{0,1,\ldots,n-1\}\).

\section{Congruence splitting lattices; uniform
refinement properties}

We shall recall in this section some of the definitions and
results of \cite{Wehr1,Wehr2} as
well as a few new ones. Recall first \cite[Section 3]{Wehr2}
that a lattice $L$ is \emph{congruence splitting} when for all
\(a\leq b\) in $L$ and all congruences $\ba_0$ and $\ba_1$ in
$L$, if \(\Theta(a,b)=\ba_0\vee\ba_1\), then there exist
elements $x_0$ and $x_1$ of \([a,\,b]\) such that
\(x_0\vee x_1=b\) and, for all \(i<2\),
\(\Theta(a,x_i)\subseteq\ba_i\).

In \cite[Proposition 3.2]{Wehr2}, we give a list of
sufficient conditions for a lattice to be congruence
splitting; this can be recorded here in the following fashion:

\begin{proposition}\label{P:SuffCS}
The following holds:
\begin{itemize}

\item[\rm (a)] Every lattice that is either relatively
complemented or sectionally complemented is congruence
splitting.

\item[\rm (b)] Every atomistic lattice is congruence
splitting.

\item[\rm (c)] The class of congruence splitting lattices is
closed under direct limits.\qed

\end{itemize}
\end{proposition}

There are easy
examples of non congruence splitting lattices, as for
example any chain with at least three elements. However, it is
to be noted that two lattices may have isomorphic
congruence lattices while one is congruence
splitting and the other one is not. Our next definition is
related to the effect of the congruence splitting
property on the congruence lattice alone.

\begin{definition}\label{D:URP}
(see \cite[Definition 2.1]{Wehr2})
Let $S$ be a semilattice, let $e$ be an element of $S$.
Say that the \emph{uniform refinement property} (URP) holds at
$e$ when for all families \((a_i)_{i\in I}\) and
\((b_i)_{i\in I}\) of elements of \(S\) such that
\(a_i\vee b_i=e\) (all \(i\in I\)), there are families
\((a^*_i)_{i\in I}\), \((b^*_i)_{i\in I}\) and 
\((c_{ij})_{(i,j)\in I\times I}\) of elements of \(S\)
such that for all \(i,j,k\in I\), we have
\begin{itemize}
\item[(i)] \(a^*_i\leq a_i\) and
\(b^*_i\leq b_i\) and \(a^*_i\vee b^*_i=e\).

\item[(ii)] \(c_{ij}\leq a^*_i,b^*_j\) and
\(a^*_i\leq a^*_j\vee c_{ij}\).

\item[(iii)] \(c_{ik}\leq c_{ij}\vee c_{jk}\).
\end{itemize}

Say that $S$ satisfies the URP when the URP holds at every
element of $S$.

Then define similarly the \emph{weak uniform refinement
property} (WURP) \emph{at $e$} when, under the same hypotheses
on the $a_i$, $b_i$ ($i\in I$) and $e$, there are \(c_{ij}\)
such that for all \(i,j,k\in I\), we have
\begin{itemize}
\item[(i$'$)] \(c_{ij}\leq a_i,b_j\).

\item[(ii$'$)] \(c_{ij}\vee a_j\vee b_i=e\).

\item[(iii$'$)] \(c_{ik}\leq c_{ij}\vee c_{jk}\).
\end{itemize}
Say that $S$ satisfies the WURP when the WURP holds at
every element of $S$.
\end{definition}

\begin{lemma}\label{L:twoURPs}
In the context above, the URP implies the WURP.\qed
\end{lemma}

\begin{proposition}\label{P:DistrURP}
Every distributive lattice satisfies the URP.
\end{proposition}

\begin{proof}
Let \((a_i)_{i\in I}\) and \((b_i)_{i\in I}\) two families of
elements of a distributive lattice $D$ such that
\(a_i+b_i=\mathrm{constant}\). It is easy to verify that
the elements \(a_i^*=a_i\), \(b_i^*=b_i\) and
\(c_{ij}=a_i\wedge b_j\) are as required.
\end{proof}

We will need later the following straightforward lemma (see
also \cite[Proposition 2.3]{Wehr2} for the URP):

\begin{lemma}\label{L:WDUnif}
Let \(\mu:\,S\to T\) be a weak-distributive homomorphism of
semilattices and let $e$ be an element of $S$. If URP
(resp. WURP) holds at $e$ in $S$, then it also holds at
\(\mu(e)\) in $T$.\qed
\end{lemma}

\begin{corollary}\label{C:WDImDistr}
Let $S$ be a distributive semilattice. If $S$ is the
image of a distributive lattice by a weak-distributive
homomorphism, then $S$ satisfies the URP (thus the WURP).\qed
\end{corollary}

In particular, any distributive semilattice that is the image
of a generalized Boolean algebra under a weak-distributive
homomorphism (this is E. T. Schmidt's sufficient
condition for being isomorphic to the congruence semilattice
of a lattice, \cite{Schm68,Tisc94}) satisfies the WURP.

We end this section by recording one of the main results of
\cite{Wehr2}:

\begin{theorem}\label{T:CSURP}
{\rm \cite[Theorem 3.3]{Wehr2}}
Let $L$ be a congruence splitting lattice. Then \(\C(L)\)
satisfies the URP (thus the WURP).\qed
\end{theorem}

Hence, by Lemma~\ref{L:twoURPs}, if $L$ is a congruence
splitting lattice, then \(\C(L)\) also has the WURP. In
particular, if we manage to find a lattice such that its
congruence semilattice does not satisfy the WURP, then this
lattice cannot be embedded con\-gru\-en\-ce-pre\-ser\-ving\-ly
into a congruence splitting lattice.

\section{The decorations of $M_3$ and $N_5$}

From now on until Theorem~\ref{T:NonUnif}, we shall fix a
non distributive lattice variety $\V$. Let $\ch$ denote the
two-element chain. For every set $X$, denote by $\E(X)$ the
free product (=coproduct) of $X$ copies of $\ch$ in $\V$.
Denote by $\B(X)$ the bounded lattice obtained from $\E(X)$
by adding a new largest element $1$ and a new least element
$0$; write \(\E_{\V}(X)\) (resp. \(\B_{\V}(X)\)) if $\V$
needs to be specified.
Thus $\B(X)$ is generated as a bounded lattice by chains
\(s_i<t_i\) (all $i\in X$).
Note that if $Y$ is a subset of $X$, then there is
a canonical retraction from $\B(X)$ onto $\B(Y)$, sending each
$s_i$ (resp. $t_i$) to $0$ for every \(i\in X\setminus Y\).
Thus, we shall often identify $\B(Y)$ with the bounded
sublattice of $\B(X)$ generated by all $s_i$ and $t_i$
($i\in Y$). Moreover, the abovementioned retraction from
$\B(X)$ onto $\B(Y)$ induces a retraction from $\C(\B(X))$
onto $\C(\B(Y))$. Hence, we shall also identify $\C(\B(Y))$
with the corresponding subsemilattice of $\C(\B(X))$.

From now on until Theorem~\ref{T:NonUnif}, we shall fix a set
$X$ such that
\(|X|\geq\aleph_2\). We denote, for all $i\in X$, by $\ba_i$
and $\bb_i$ the compact congruences of $\B(X)$ defined by
\begin{equation}
\ba_i=\Theta(0,s_i)\vee\Theta(t_i,1);\qquad
\bb_i=\Theta(s_i,t_i).
\end{equation}
In particular, note that \(\ba_i\vee\bb_i=\one\).
Now, towards a contradiction, suppose that there are compact
congruences $\bc_{ij}$ ($i,j\in X$) such that for all
\(i,j,k\in X\), the following holds:
\begin{gather}
\bc_{ij}\subseteq\ba_i,\bb_j\label{Eq:E1}\\
\bc_{ij}\vee\ba_j\vee\bb_i=\one\label{Eq:E2}\\
\bc_{ik}\subseteq\bc_{ij}\vee\bc_{jk}.\label{Eq:E3}
\end{gather}

Since the $\C$ functor preserves direct limits, there exists,
for all $i,j\in X$, a finite subset $U=F(\{i,j\})$ of $X$ such
that both $\bc_{ij}$ and $\bc_{ji}$ belong to $\C(\B(U))$.
By Kuratowski's Theorem, there are mutually
distinct elements, which we may denote by $0$, $1$, $2$ of
$X$ such that \(0\notin F(\{1,2\})\), \(1\notin F(\{0,2\})\),
and \(2\notin F(\{0,1\})\). Let
\(\pi:\,\B(X)\twoheadrightarrow\B(3)\) be the canonical
retraction. For every $i<3$, denote by $i'$ and
$i''$ the other two elements of $3$, arranged in such a way
that \(i'<i''\). For all $i<3$, put
\(\bd_i=\C(\pi)(\bc_{i'i''})\).

Therefore, applying the semilattice homomorphism \(\C(\pi)\)
to the inequalities (\ref{Eq:E1} -- \ref{Eq:E3}) yields
\begin{gather}
\bd_0\subseteq \ba_1,\bb_2;\qquad
\bd_1\subseteq\ba_0,\bb_2;\qquad
\bd_2\subseteq\ba_0,\bb_1;\label{Eq:D1}\\
\bd_0\vee\ba_2\vee\bb_1=
\bd_1\vee\ba_2\vee\bb_0=
\bd_2\vee\ba_1\vee\bb_0=\one;\label{Eq:D2}\\
\bd_1\subseteq\bd_0\vee\bd_2.\label{Eq:D3}
\end{gather}

\begin{lemma}\label{L:Supp}
For all $i<3$, $\bd_i$ belongs to $\C(\B(3\setminus\{i\}))$.
\end{lemma}

\begin{proof}
For example for $i=0$. Since \(0\notin F(\{1,2\})\),
$\bc_{12}$ belongs to \(\B(X\setminus\{0\})\), hence
\(\bd_0\in\B(\{1,2\})\).
\end{proof}

Now, since $\V$ is a non distributive variety of lattices, by
a classical result of lattice theory, either the diamond $M_3$
or the pentagon $N_5$ belongs to the variety $\V$. Denote by
$M$ one of these lattices that belongs to $\V$ and decorate it
with three $2$-element chains \(x_i<y_i\) ($i<3$) in the
following way:\medskip

\noindent{\bf Case 1.} $M=M_3$. Let $p,q,r$ be the three atoms
of $M_3$. Put
\begin{gather*}
x_0=0\qquad y_0=p\\
x_1=q\qquad y_1=1\\
x_2=0\qquad y_2=r
\end{gather*}

\noindent{\bf Case 2.} $M=N_5$. Let $a>c$ and $b$ be the three
join-irreducible elements of $N_5$. Put
\begin{gather*}
x_0=0\qquad y_0=c\\
x_1=b\qquad y_1=1\\
x_2=0\qquad y_2=a
\end{gather*}

Both cases can be described by the following picture:

\begin{picture}(200,180)(-80,-110)
\thicklines
\put(-47.88,2.12){\line(1,1){45.76}}
\put(-47.88,-2.12){\line(1,-1){45.76}}
\put(0,3){\line(0,1){44}}
\put(0,-47){\line(0,1){44}}
\put(2.12,47.88){\line(1,-1){45.76}}
\put(2.12,-47.88){\line(1,1){45.76}}

\put(-50,0){\circle{6}}
\put(0,0){\circle{6}}
\put(50,0){\circle{6}}
\put(0,50){\circle{6}}
\put(0,-50){\circle{6}}

\put(-55,0){\makebox(0,0)[r]{$y_0=p$}}
\put(-5,0){\makebox(0,0)[r]{$x_1=q$}}
\put(55,0){\makebox(0,0)[l]{$y_2=r$}}
\put(0,-60){\makebox(0,0){$x_0=x_2=0$}}
\put(0,60){\makebox(0,0){$y_1=1$}}
\put(0,-90){\makebox(0,0){{\bf Case 1.} $M=M_3$}}

\put(177.88,-47.88){\line(-1,1){25.76}}
\put(150,-17){\line(0,1){34}}
\put(152.12,22.12){\line(1,1){25.76}}
\put(182.12,47.88){\line(1,-1){45.76}}
\put(227.88,-2.12){\line(-1,-1){45.76}}

\put(180,-50){\circle{6}}
\put(150,-20){\circle{6}}
\put(150,20){\circle{6}}
\put(180,50){\circle{6}}
\put(230,0){\circle{6}}

\put(145,-20){\makebox(0,0)[r]{$y_0=c$}}
\put(180,-60){\makebox(0,0){$x_0=x_2=0$}}
\put(145,20){\makebox(0,0)[r]{$y_2=a$}}
\put(180,60){\makebox(0,0){$y_1=1$}}
\put(235,0){\makebox(0,0)[l]{$x_1=b$}}
\put(180,-90){\makebox(0,0){{\bf Case 2.} $M=N_5$}}

\end{picture}

The relevant properties of these decorations are summarized in
both following straightforward lemmas:

\begin{lemma}\label{L:Ineq}
The decorations defined above satisfy the following
inequalities
\[
x_0\wedge y_1\leq x_1;\qquad y_1\leq x_1\vee y_0;\qquad
x_1\wedge y_2\leq x_2;\qquad y_2\leq x_2\vee y_1,
\]
but \(y_2\not\leq x_2\vee y_0\).\qed
\end{lemma}

\begin{lemma}\label{L:ChDistr}
For all $i<3$, the sublattice of $M$ generated by the elements
\(x_j\) and \(y_j\) ($j\ne i$) is distributive.\qed
\end{lemma}

This, along with (\ref{Eq:D1} -- \ref{Eq:D3}), will be
sufficient to obtain a contradiction. Note already that since
the free product of three $2$-element chains in the variety
generated by either $M_3$ or $N_5$ is finite, the problem is
already reduced to a ``computable" level. However, the size of
the corresponding computations is such that it is useful to
reduce (greatly) their complexity to mere computations in
$M_3$ and $N_5$. This is what we shall do in
Section~\ref{Dlat}.

\section{Reduction to the distributive world}\label{Dlat}

From now on, we shall denote by $D$ be the free product of
two $2$-element chains in the variety of all
\emph{distributive} lattices. Thus, $D$ is generated by two
chains \(u_0<v_0\) and \(u_1<v_1\).

Thus $D$ is a finite distributive lattice, that can be
represented by the following diagram:

\begin{picture}(200,300)(-100,-170)
\thicklines
\put(0,45){\circle{6}}
\put(30,15){\circle{6}}
\put(30,75){\circle{6}}
\put(60,15){\circle{6}}
\put(60,45){\circle{6}}
\put(60,105){\circle{6}}
\put(90,15){\circle{6}}
\put(90,75){\circle{6}}
\put(120,45){\circle{6}}

\put(0,-45){\circle{6}}
\put(30,-15){\circle{6}}
\put(30,-75){\circle{6}}
\put(60,-15){\circle{6}}
\put(60,-45){\circle{6}}
\put(60,-105){\circle{6}}
\put(90,-15){\circle{6}}
\put(90,-75){\circle{6}}
\put(120,-45){\circle{6}}

\put(2.12,47.12){\line(1,1){25.76}}
\put(32.12,77.12){\line(1,1){25.76}}
\put(32.12,72.88){\line(1,-1){25.76}}
\put(2.12,42.12){\line(1,-1){25.76}}
\put(32.12,17.12){\line(1,1){25.76}}
\put(62.12,47.12){\line(1,1){25.76}}
\put(92.12,72.88){\line(1,-1){25.76}}
\put(62.12,42.88){\line(1,-1){25.76}}
\put(92.12,17.12){\line(1,1){25.76}}
\put(62.12,102.88){\line(1,-1){25.76}}

\put(2.12,-47.12){\line(1,-1){25.76}}
\put(32.12,-77.12){\line(1,-1){25.76}}
\put(32.12,-72.88){\line(1,1){25.76}}
\put(2.12,-42.12){\line(1,1){25.76}}
\put(32.12,-17.12){\line(1,-1){25.76}}
\put(62.12,-47.12){\line(1,-1){25.76}}
\put(92.12,-72.88){\line(1,1){25.76}}
\put(62.12,-42.88){\line(1,1){25.76}}
\put(92.12,-17.12){\line(1,-1){25.76}}
\put(62.12,-102.88){\line(1,1){25.76}}

\put(32.12,12.88){\line(1,-1){25.76}}
\put(62.12,12.88){\line(1,-1){25.76}}
\put(32.12,-12.88){\line(1,1){25.76}}
\put(62.12,-12.88){\line(1,1){25.76}}
\put(30,-12){\line(0,1){24}}
\put(90,-12){\line(0,1){24}}
\put(60,18){\line(0,1){24}}
\put(60,-42){\line(0,1){24}}

\put(-6,-45){\makebox(0,0)[r]{$u_0$}}
\put(-6,45){\makebox(0,0)[r]{$v_0$}}
\put(126,-45){\makebox(0,0)[l]{$u_1$}}
\put(126,45){\makebox(0,0)[l]{$v_1$}}
\put(60,-115){\makebox(0,0){$u_0\wedge u_1=0$}}
\put(60,115){\makebox(0,0){$v_0\vee v_1=1$}}

\put(80,-140){\makebox(0,0){{\bf The lattice}\ $D$}}

\end{picture}

For all \(i<3\), let \(\pi_i:\,\B(3\setminus\{i\})\to D\) be
the unique lattice homomorphism sending
\(s_{i'}\) to \(u_0\), \(t_{i'}\) to \(v_0\),
\(s_{i''}\) to \(u_1\), \(t_{i''}\) to \(v_1\). Furthermore,
let \(\rho:\,\B(3)\to M\) be the unique lattice
homomorphism sending $s_i$ to $x_i$ and $t_i$ to $y_i$ (all
$i<3$); let
$\rho_i$ be the restriction of $\rho$ to
\(\B(3\setminus\{i\})\).

\begin{lemma}\label{L:BasicDistr}
Let $L$ be any distributive lattice, let $a$, $b$, $a'$, $b'$
be elements of $L$. Then we have
\[
\Theta^+(a,b)\cap\Theta^+(a',b')=
\Theta^+(a\wedge a',b\vee b').
\]
\end{lemma}

\begin{proof}
Let $B$ be the generalized Boolean algebra R-generated by $L$
(in the sense of \cite[Part II, Section 4]{Grat78}); identify
every congruence $\theta$ on $L$ with the unique congruence on
$B$ extending $\theta$. For all elements $x$ and $y$ of $B$,
denote by \(x\sd y\) the unique relative complement of
\(x\wedge y\) in the interval \([0,\,x]\), and then put
\(x\triangle y=(x\sd y)\vee(y\sd x)\).
Then a pair $(x,y)$
belongs to \(\Theta^+(a,b)\) (resp. \(\Theta^+(a',b')\)) if
and only if \(x\triangle y\leq a\sd b\) (resp.
\(x\triangle y\leq a'\sd b'\)). Therefore, $(x,y)$ belongs to
\(\Theta^+(a,b)\cap\Theta^+(a',b')\) if and only if
\(x\triangle y\leq(a\sd b)\wedge(a'\sd b')=
(a\wedge a')\sd(b\vee b')\).
\end{proof}

\begin{remark}
In particular, one recovers the classical result
that if \(a\leq b\leq c\leq d\) are elements
of any distributive lattice, then
\(\Theta(a,b)\cap\Theta(c,d)=\zero\).
\end{remark}

Now, for all $i<3$, put \(\be_i=\C(\pi_i)(\bd_i)\).

\begin{lemma}\label{L:eiUnique}
For all \(i<3\),
\(\be_i=
\Theta^+(u_0\wedge v_1,u_1)\vee\Theta^+(v_1,u_1\vee v_0)\).
\end{lemma}

\begin{proof}
Applying $\C(\pi_i)$ to
the inequalities (\ref{Eq:D1}) and (\ref{Eq:D2}) yields
both following inequalities:
\begin{gather}
\be_i\subseteq\Theta(0,u_0)\vee\Theta(v_0,1)\quad
{\rm and}\quad\be_i\subseteq\Theta(u_1,v_1)\label{Eq:euv1}\\
\be_i\vee\Theta(0,u_1)\vee\Theta(v_1,1)\vee\Theta(u_0,v_0)=
\one.\label{Eq:euv2}
\end{gather}
However, $D$ is a finite distributive lattice, thus $\C(D)$ is
a finite Boolean algebra, and, by the Remark above, for all
\(j<2\), the elements \(\Theta(0,u_j)\vee\Theta(v_j,1)\)
and \(\Theta(u_j,v_j)\) are complemented elements of
$\C(D)$; in fact, $\one$ is the disjoint union of
\(\Theta(0,u_j)\), \(\Theta(u_j,v_j)\), \(\Theta(v_j,1)\).
Then, from both inequalities (\ref{Eq:euv1}) and
(\ref{Eq:euv2}) and a new application of
Lemma~\ref{L:BasicDistr}, one deduces easily that
\begin{align*}
\be_i&=(\Theta(0,u_0)\vee\Theta(v_0,1))\cap\Theta(u_1,v_1)\\
&=(\Theta(0,u_0)\cap\Theta(u_1,v_1))\vee
(\Theta(v_0,1)\cap\Theta(u_1,v_1))\\
&=\Theta^+(u_0\wedge v_1,u_1)\vee
\Theta^+(v_1,u_1\vee v_0).\tag*{\qed}
\end{align*}
\renewcommand{\qed}{}
\end{proof}

Now, for all $i<3$, it results from Lemma~\ref{L:ChDistr} that
there exists a unique lattice homomorphism
\(\varphi_i:\,D\to M\) such that
\(\varphi_i\circ\pi_i=\rho_i\). The corresponding commutative
diagram is the following:

\begin{picture}(100,100)(-130,-20)

\put(20,10){\vector(2,1){80}}
\put(120,50){\vector(0,-1){40}}
\put(50,0){\vector(1,0){40}}

\put(60,40){\makebox(0,0)[tr]{$\pi_i$}}
\put(125,30){\makebox(0,0)[l]{$\varphi_i$}}
\put(70,-5){\makebox(0,0)[t]{$\rho_i$}}

\put(-20,0){\makebox(0,0)
{$(\B(3\setminus\{i\}),s_{i'},t_{i'},s_{i''},t_{i''})$}}
\put(140,0){\makebox(0,0)
{$(M,x_{i'},y_{i'},x_{i''},y_{i''})$}}
\put(140,60){\makebox(0,0)
{$(D,u_0,v_0,u_1,v_1)$}}

\end{picture}

Since $\C$ is a functor, we get from this and from
Lemma~\ref{L:eiUnique} that for all $i<3$, we have
\begin{multline}
\C(\rho)(\bd_i)=\C(\varphi_i)(\be_i)=
\C(\varphi_i)
(\Theta^+(u_0\wedge v_1,u_1)\vee\Theta^+(v_1,u_1\vee v_0))\\
=\Theta^+(x_{i'}\wedge y_{i''},x_{i''})\vee
\Theta^+(y_{i''},x_{i''}\vee y_{i'}).
\end{multline}
In particular, we have, using Lemma~\ref{L:Ineq},
\begin{align*}
\C(\rho)(\bd_0)&=
\Theta^+(x_1\wedge y_2,x_2)\vee\Theta^+(y_2,x_2\vee y_1)
=\zero,\\
\C(\rho)(\bd_2)&=
\Theta^+(x_0\wedge y_1,x_1)\vee\Theta^+(y_1,x_1\vee y_0)
=\zero,\\
{\rm but}\qquad\C(\rho)(\bd_1)&=
\Theta^+(x_0\wedge y_2,x_2)\vee\Theta^+(y_2,x_2\vee y_0)
\ne\zero.\\
\end{align*}
On the other hand, by applying $\C(\rho)$ to (\ref{Eq:D3}), we
obtain that
\[
\C(\rho)(\bd_1)\leq\C(\rho)(\bd_0)\vee\C(\rho)(\bd_2),
\]
a contradiction. Therefore, we have proved the following
theorem:

\begin{theorem}\label{T:NonUnif}
Let $\V$ be any non distributive variety of lattices, let $X$
be any set such that $|X|\geq\aleph_2$. Let $\B_{\V}(X)$ be the
free product in $\V$ of $X$ copies of a $2$-element chain with
a least and a largest element added. Then \(\C(\B_{\V}(X))\)
does not satisfy WURP at $\one$.\qed
\end{theorem}

\section{Extensions to further lattices and to regular rings}

This section will be devoted to harvest consequences of
Theorem~\ref{T:NonUnif}.

\begin{corollary}\label{C:FreeNon}
Let $L$ be any lattice that admits a lattice homomorphism
onto a free bounded lattice in the variety generated by either
\(M_3\) or \(N_5\) with \(\aleph_2\) generators.
Then \(\C(L)\) does not satisfy WURP at
$\one$. In particular, there exists no congruence
splitting lattice $K$ such that \(\C(K)\cong\C(L)\);
furthermore, \(\C(L)\) does not satisfy Schmidt's condition.
\end{corollary}

\begin{proof}
Let \(\V\) be the lattice variety generated by either
\(M_3\) or \(N_5\) and, for any set $X$, let \(\F_{\V}(X)\) be
the free bounded lattice on $X$ in $\V$. First, if the
cardinality of $X$ is infinite, then there exists
a surjective lattice homomorphism from \(\F_{\V}(X)\) onto
\(\B_{\V}(X)\) (split $X$ into two disjoint sets $X_0$ and
$X_1$ such that \(|X_0|=|X_1|=|X|\); send the elements of
$X_1$ (resp. $X_2$) onto all $s_i$ (resp. $t_i$)). Therefore,
if \(|X|=\aleph_2\), then there exists by assumption a
surjective lattice homomorphism
\(f:\,L\twoheadrightarrow\B_{\V}(X)\). By
\cite[Proposition 1.2]{Wehr2}, the corresponding congruence
mapping \(\C(f):\,\C(L)\twoheadrightarrow\C(\B_{\V}(X))\) is
weak-distributive. Thus, if $\C(L)$ would satisfy the WURP at
$\one$, then, by Lemma~\ref{L:WDUnif}, $\C(\B_{\V}(X))$ would
also satisfy the same refinement property, therefore
contradicting Theorem~\ref{T:NonUnif}. The last two statements
result from Theorem~\ref{T:CSURP} and
Corollary~\ref{C:WDImDistr}.
\end{proof}

This shows in particular that there are distributive
semilattices that are representable as congruence
semilattices of lattices (the \(\C(L)\)'s, with $L$ free
lattice on at least $\aleph_2$ generators in any
non-distributive variety) that, nevertheless, do not satisfy
any of the known sufficient conditions implying
representability (as Schmidt's condition).

\begin{corollary}\label{C:FreeNon2}
Let $\V$ be any non distributive variety of lattices
and let $F$ be any free (resp. free bounded) lattice with at
least $\aleph_2$ generators in $\V$. Then there exists no
congruence splitting lattice $K$ such that
\(\C(K)\cong\C(F)\).\qed
\end{corollary}

\begin{corollary}\label{C:RegRings}
Let $\V$ and $F$ as above. Then there exists no von Neumann
regular ring $R$ whose semilattice of finitely generated
two-sided ideals is isomorphic to \(\C(F)\).
\end{corollary}

\begin{proof}
Proposition~\ref{P:SuffCS} (and the fact that the lattice of
principal right ideals of $R$ is sectionally complemented),
Corollary~\ref{C:FreeNon2} and \cite[Corollary 4.4]{Wehr2}.
\end{proof}

We do not know whether every lattice of cardinality at most
$\aleph_1$ admits a con\-gru\-en\-ce-pre\-ser\-ving extension
into a sectionally complemented lattice. On the other hand, G.
Gr\"atzer and E. T. Schmidt prove in \cite{GrSc} that every
\emph{finite} lattice has a con\-gru\-en\-ce-pre\-ser\-ving
extension into a sectionally complemented finite lattice. Note
that if
\(\V\) is a non distributive lattice variety generated by a
single finite lattice, then \(\F_{\V}(X)\) is a direct limit
of a limit system of finite lattices with embeddings
having the congruence extension property; nevertheless, its
congruence semilattice is complicated in the sense that it
does not satisfy WURP.

\section*{Acknowledgments}

The authors wish to
thank warmly Ralph Freese who, when questioned about it,
returned to them very quickly the construction of the lattice
freely generated by two $2$-element chains in the pentagon
variety. Although its use has eventually been eliminated in
the final version of the proof, it greatly helped in its
elaboration. The paper was completed while 
the second author was visiting the University of Caen. The
hospitality and excellent conditions provided by the
mathematics department are greatly appreciated.

\end{document}